\documentclass{article}[12pt]
\usepackage{amssymb}
\usepackage{eurosym}
\usepackage{epsfig}
\usepackage{amsmath}
\usepackage{amsfonts}
\usepackage{pgf,tikz}
\usepackage{enumerate}
\usepackage{cite}
\usepackage[colorlinks=true, linkcolor=black, citecolor=black, urlcolor=black]{hyperref}

\newcommand{\R}{{{\Bbb R}}}

\newtheorem{theorem}{\sc Theorem}[section]
\newtheorem{proposition}[theorem]{\sc Proposition}
\newtheorem{definition}[theorem]{\sc Definition}
\newtheorem{remark}[theorem]{\sc Remark}
\newtheorem{corollary}[theorem]{\sc Corollary}
\newtheorem{example}[theorem]{\sc Example}

\def\qed{\hbox to 0pt{}\hfill$\rlap{$\sqcap$}\sqcup$\medbreak}

\title{A Birkhoff--Kellogg type theorem for discontinuous operators with applications}

\author{Alessandro Calamai,
Gennaro Infante and Jorge Rodr\'iguez--L\'opez}

\date{}
\begin{document}
 \maketitle


\medbreak

\noindent {\it Abstract.} 
By means of fixed point index theory for multivalued maps, 
we provide an analogue of the classical Birkhoff--Kellogg Theorem in the context of discontinuous operators acting on affine wedges in Banach spaces. Our theory is fairly general and can be applied, for example, to eigenvalues and parameter problems for ordinary differential equations with discontinuities. We illustrate in details this fact for a class of second order boundary value problem with deviated arguments and discontinuous terms. In a specific example, we explicitly compute the terms that occur in our theory. 

\medbreak


\medbreak


\section{Introduction}

The celebrated invariant-direction Theorem due to Birkhoff and Kellogg~\cite{B-K-1922} is an abstract existence result
 that, roughly speaking, gives conditions for the existence of a ``nonlinear'' eigenvalue and eigenvector for compact maps in normed linear spaces.
Among its various extensions, one is set in cones and is due to Krasnosel'ski\u{i} and Lady\v{z}enski\u{\i}~\cite{Kra-Lady}.
These classical functional analytic tools find applications e.g.\ to eigenvalue problems for ODEs and PDEs (see for example the book~\cite{ADV}, the recent papers~\cite{gi-BK, giaml22} and references therein); typically, the methodology in this context is to reformulate the  
given boundary value problem as a fixed point problem in a suitable Banach space.

Recently the first two authors developed a Birkhoff--Kellogg type theorem in the framework of
\textit{affine cones} (cf.\ \cite{CI_MMAS}, see also \cite{acgi16, acgi3, djeb2014}). The motivation for this new type of results is that the setting of affine cones seems to be helpful when dealing with equations with delay effects. A key ingredient in \cite{CI_MMAS} is the \emph{continuity} of the involved operator. On the other hand, there has been recently a rising attention towards \textit{discontinuous} differential equations, that occur when modelling real world phenomena. Here we mention the classical books by Filippov \cite{Filippov}, Carl and Heikkil\"a \cite{CarlHeik}, and Heikkil\"a  and  Lakshmikantham \cite{HeLak} and the more recent book by Figueroa, Pouso and RL\cite{FLR}.

In the present paper we provide a discontinuous version of the Birkhoff--Kellogg type result
in the setting of affine wedges in Banach spaces,
see Theorem \ref{th2_BK} below.
The proof of Theorem \ref{th2_BK} is based on
 the fixed point index theory for discontinuous operators developed in \cite{FLR}.
We stress that a crucial point in the construction of the index for discontinuous operators is its equivalence with the corresponding one of a suitable \textit{multivalued} map, for which it is already defined, see \cite{fp1}. Note that this newly constructed topological tool for discontinuous operators inherits the key properties of the classical one. This construction is sketched in Section \ref{sec-BK} for completeness.

In Section \ref{sec-appl} we illustrate the applicability of our results to boundary value problems,
see Theorem~\ref{th}.  In more details, we consider the following second order parameter-dependent differential equation with deviated argument
\begin{equation}\label{eq_2or-intro}
	u''(t)+\lambda\,f(t,u(t),u(\sigma(t)))=0, \quad t\in[0,1],
\end{equation}
with initial condition
\begin{equation}\label{eq_initial-intro}
	u(t)=\omega(t), \quad t\in[-r,0],
\end{equation}
and the final homogeneous boundary condition
\begin{equation}\label{BC-intro}
	u(1)=0,
\end{equation}
where $\lambda \geq 0$ is a parameter, $r\geq 0$, $\sigma:[0,1]\rightarrow[-r,1]$ and $\omega:[-r,0]\rightarrow[0,\infty)$ are suitable continuous functions, while the nonlinearity $f:[0,1]\times[0,\infty)\times[0,\infty)\rightarrow[0,\infty)$ may be discontinuous with respect to the second argument in an appropriate sense.
We employ a concept of \textit{admissible discontinuity curve} as in \cite{FLR}.
We conclude the paper by illustrating the applicability of our theory by means of a toy model with delay, see Example \ref{example2}.

As far as we know, our results extend and complement the previous literature. This is highlighted in more details in Remarks \ref{rem1}, \ref{rem2} and \ref{eigen-clas}.

\section{Birkhoff-Kellogg type results via fixed point index theory}\label{sec-BK}

\subsection{On fixed point index theory for discontinuous operators}

Let $K$ be a nonempty closed and convex subset of a real Banach space $(X,\| \cdot \|)$, $U\subset K$ a relatively open subset and $T:\overline{U}\subset K \longrightarrow K$ a mapping, not necessarily continuous.  
\begin{definition}
	\label{def1}
	The closed--convex envelope of an operator $T:\overline{U}\subset K \longrightarrow K$ is the multivalued mapping  $\mathbb{T}: \overline{U} \longrightarrow 2^X$ given by
	\begin{equation}\label{TT}
		\mathbb{T}x=\bigcap_{\varepsilon>0}\overline{\rm co} \, T\left(\overline{B}_{\varepsilon}(x)\cap \overline{U}\right) \quad \text{for every } x\in \overline{U},
	\end{equation}
	where $\overline{B}_{\varepsilon}(x)$ denotes the closed ball centered at $x$ and radius $\varepsilon$, and $\overline{\rm co}$ means closed convex hull.  
\end{definition}	

\begin{example}
	\begin{enumerate}
		\item Consider the real function $T:\mathbb{R}\to\mathbb{R}$ defined as $T(x)=x$, if $x\leq 0$, and $T(x)=x+1$, if $x>0$. Its closed-convex envelope is the multivalued map $\mathbb{T}$ given by $\mathbb{T}(x)=\{x\}$, if $x<0$; $\mathbb{T}(x)=\{x+1\}$, if $x>0$; and $\mathbb{T}(0)=[0,1]$.
		\item The closed-convex envelope of any continuous map $T$ is equal to $T$. 
	\end{enumerate}
\end{example}

Now we recall some useful properties of closed--convex envelopes (cc--envelopes for short) and the definition of the fixed point index that we will employ throughout this paper. The reader is referred to \cite{DegDisc,FLR} for details.
\begin{proposition}
	\label{pro2}
	Let $\mathbb T$ be the cc--envelope of an operator $T: \overline{U} \longrightarrow K$. Then the following properties hold:
	
	\begin{enumerate}
		
		\item If  $\tilde {\mathbb T}:\overline{U} \longrightarrow 2^X$ is an upper semicontinuous (usc) operator which assumes closed and convex values and $Tx \in \tilde{\mathbb T} x$ for all $x \in \overline{U}$, then $\mathbb T x \subset \tilde {\mathbb T} x$ for all $x \in \overline{U}$;
		
		\item If $T$ maps bounded sets into relatively compact sets, then $\mathbb T$ assumes compact values and it is usc;
		\item If $T \, \overline{U}$ is relatively compact, then $\mathbb T \, \overline{U}$ is relatively compact.
		
	\end{enumerate}
	
\end{proposition}

The fixed point index for a not necessarily continuous operator $T$ was introduced in \cite{FLR_JFPTA} by using the degree theory developed in \cite{DegDisc} and a retraction trick, just as in the classical case. Both topological degree and fixed point index theories are based on the available results for the multivalued cc-envelope $\mathbb{T}$.

\begin{definition}\label{def_index}
	Let $T:\overline{U}\subset K \longrightarrow K$ be an operator such that $T\,\overline{U}$ is relatively compact, $T$ has no fixed points on $\partial\, U$ and 
	\begin{equation}\label{cond}
		\left\{x\right\}\cap\mathbb{T}x\subset\left\{Tx\right\} \quad \mbox{for every $x\in\overline{U} \cap \mathbb{T}\,\overline{U}$},
	\end{equation}
	where $\mathbb T$ is the cc--envelope of $T$.	
	
	We define the fixed point index of $T$ in $K$ over $U$ as
	\[i_{K}(T,U)=\deg(I-T\circ r,r^{-1}(U),0),\]
	where $r$ is a continuous retraction of $X$ onto $K$ and $\deg$ is the degree introduced in  \cite{DegDisc}. 
\end{definition}	

\begin{remark}
	Note that condition \eqref{cond} means that the set of fixed points of $\mathbb{T}$ (i.e., the set of points $x$ such that $x\in \mathbb{T}x$) is contained in the set of fixed points of $T$. This is a weaker condition than the continuity of $T$; indeed, if $T$ is continuous, then $\mathbb{T}x=\{Tx\}$ for all $x\in \overline{U}$ and thus \eqref{cond} is trivially satisfied.
\end{remark}

We now recall a useful proposition from \cite{FLR_JFPTA} that relates 
the fixed point index of the discontinuous operator $T$ with that of its associated multivalued mapping $\mathbb{T}$.

\begin{proposition}\cite[Proposition 2.12]{FLR_JFPTA}\label{prop_indT=indTT}
	Let $T$ be a mapping that satisfies the conditions of \hyperref[def_index]{Definition \ref*{def_index}}. Then, the fixed point index of $T$ is such that \[i_{K}(T,U)=i_{K}(\mathbb{T},U),\]
	where the right-hand index is the fixed point index defined for multivalued mappings, see \cite{fp1}.
\end{proposition}

As a straightforward consequence of the fixed point index theory for usc multivalued mappings, the following properties can be derived (see \cite{FLR}).

\begin{theorem}\label{th_prop_index}
	Let $T$ be a mapping that satisfies the conditions of \hyperref[def_index]{Definition \ref*{def_index}}. Then the following properties hold: 
	\begin{enumerate}[i.]
		\item (Homotopy invariance) Let $H:\overline{U}\times\left[0,1\right]\longrightarrow K$ be a mapping such that:
		\begin{enumerate}[(a)]
			\item for each $(x,t)\in\overline{U}\times[0,1]$ and all $\varepsilon>0$ there exists $\delta=\delta(\varepsilon,x,t)>0$ such that 
			\[s\in[0,1], \ \left|t-s\right|<\delta \ \Longrightarrow \ \left\|H(z,t)-H(z,s)\right\|<\varepsilon \ \ \forall\,z\in\overline{B}_{\delta}(x)\cap\overline{U};\]
			\item $H\left(\overline{U}\times\left[0,1\right]\right)$ is relatively compact;
			\item $\left\{x\right\}\cap\mathbb{H}_{t}(x)\subset\left\{H_{t}(x)\right\}$ for all $t\in\left[0,1\right]$ and all $x\in\overline{U}\cap \mathbb H_t \overline{U}$, where $H_{t}(\cdot) := H(\cdot,t)$ and $\mathbb H_{t}$ denotes the cc--envelope of $H_{t}$.
		\end{enumerate}
		If $x\neq H(x,t)$ for all $(x,t)\in\partial\,U\times\left[0,1\right]$ then the index $i_{K}\left(H_{t},U\right)$ does not depend on $t\in[0,1]$.
		\item (Additivity) Let $U$ the disjoint union of two open sets $U_{1}$ and $U_{2}$. If $0\not\in\left(I-T\right)\left(\overline{U}\setminus(U_{1}\cup U_{2})\right)$, then  \[i_{K}\left(T,U\right)=i_{K}\left(T,U_{1}\right)+i_{K}\left(T,U_{2}\right).\]
		\item (Excision) Let $A\subset U$ be a closed set. If $0\not\in\left(I-T\right)\left(\partial\,U\right)\cup\left(I-T\right)(A)$, then
		\[i_{K}\left(T,U\right)=i_{K}\left(T,U\setminus A\right).\]
		\item (Existence) If $i_{K}\left(T,U\right)\neq 0$, then there exists $x\in U$ such that $Tx=x$.
		\item (Normalization) For every constant map $T$ such that $T\,\overline{U}\subset U$, $i_{K}\left(T,U\right)=1$.
	\end{enumerate}
\end{theorem}

\newpage

\subsection{Birkhoff-Kellogg theorem and discontinuous operators}

The following notions will be used along the text. A closed convex subset $K$ of a Banach space $(X,\left\|\cdot\right\|)$ is a \textit{wedge} if $\mu\,x\in K$ for every $x\in K$ and for all $\mu\geq 0$. Furthermore, if a wedge $K$ satisfies that $K\cap (-K)=\{0\}$, then it is said to be a \textit{cone}.	
A cone $K$ induces the partial order in $X$ given by $u\preceq v$ if and only if $v-u\in K$. The cone $K$ is called \textit{normal} if there
exists $c>0$ such that $\left\|u\right\|\leq c\left\|v\right\|$ for all $u,v\in X$ with $0 \preceq u \preceq v$.	  

Let $K$ be a wedge of a Banach space $(X,\left\|\cdot\right\|)$. For a given $y\in X$, the \textit{translate} of the wedge $K$ is defined as follows
\[K_y:=y+K=\{y+x:x\in K \}. \]
Given an open bounded subset $D\subset X$ with $0\in D$, we will denote $D_{K_y}:=(y+D)\cap K_y$, which is a relatively open subset of $K_y$. 
By $\overline{D}_{K_y}$ and $\partial\,D_{K_y}$ we will mean, respectively,
the closure and  the boundary of $D_{K_y}$ relative to $K_y$.

	For the convenience of the reader, we recall here the classical Birkhoff--Kellogg Theorem \cite{B-K-1922} and a variant of it set in \emph{cones}. The latter result is due to Krasnosel'ski\u{i} and Lady\v{z}enski\u{\i}~\cite{Kra-Lady} (see also \cite[Theorem 2.3.6]{guolak}).

	\begin{theorem}[Birkhoff-Kellogg]
		Let $U$ be a bounded open neighborhood of $0$ in an infinite-dimensional normed linear space $X$, and $T:\partial\,U\longrightarrow X$ a compact map satisfying $\left\|Tx\right\|\geq \alpha>0$ for all $x\in\partial\,U$. Then there exist $x_0\in\partial\,U$ and $\lambda_0>0$ such that $x_0=\lambda_0\,T x_0$.   
	\end{theorem}
	
	\begin{theorem}[Krasnosel'ski\u{i}-Lady\v{z}enski\u{\i}]\label{th_KL}
		Let $X$ be a real Banach space, $U\subset X$ be an open bounded set with $0\in U$, $K\subset X$ be a cone, $T:K\cap\overline{U}\longrightarrow K$ be compact and suppose that
		\[\inf_{x\in K\cap\partial\,U}\left\|Tx\right\|>0. \]
		Then there exist $x_0\in K\cap\partial\,U$ and $\lambda_0>0$ such that $x_0=\lambda_0\, Tx_0$.
	\end{theorem}	

In the context of \emph{affine cones} a Birkhoff-Kellogg type result was recently proved in \cite[Theorem 2]{CI_MMAS}. It reads as follows.

\begin{theorem}\label{th_BK_affine}
	Let $(X,\left\|\cdot\right\|)$ be a real Banach space, $K\subset X$ be a cone and $D\subset X$ be an open bounded set with $y\in D_{K_y}$. Assume that $T:\overline{D}_{K_y}\longrightarrow K$ is a compact map and consider the operator
	\[T_{(y,\lambda)}:=y+\lambda\,T, \qquad (\lambda\in \mathbb{R}). \]
	Assume that there exists $\bar{\lambda}>0$ such that $i_{K_y}(T_{(y,\bar{\lambda})},D_{K_y})=0$. Then there exist $x^*\in\partial\,D_{K_y}$ and $\lambda^*\in(0,\bar{\lambda})$ such that $x^*=y+\lambda^*\, T(x^*)$. 
\end{theorem}

Now we present a discontinuous version of this Birkhoff-Kellogg type result in \emph{affine wedges}.

\begin{theorem}\label{th1_BK}
	Let $D\subset X$ be an open bounded set with $0\in D$, $y\in X$ be fixed
	and $K$ be a wedge. Assume that $T:\overline{D}_{K_y} \longrightarrow K$ is a mapping such that $T\,\overline{D}_{K_y}$ is relatively compact and consider the operator
	\[T_{(y,\lambda)}:=y+\lambda\,T, \qquad (\lambda\in \mathbb{R}). \]
	Moreover, assume that there exists $\bar{\lambda}>0$ such that $i_{K_y}(T_{(y,\bar{\lambda})},D_{K_y})\neq 1$ and for each $\lambda\in(0,\bar{\lambda}]$,
	\begin{equation}\label{cond_Ty}
		\left\{x\right\}\cap\mathbb{T}_{(y,\lambda)}(x)\subset\left\{T_{(y,\lambda)}(x)\right\} \quad \mbox{for every $x\in\overline{D}_{K_y}$},
	\end{equation}
	where $\mathbb{T}_{(y,\lambda)}$ denotes the cc--envelope of $T_{(y,\lambda)}$.
	
	Then there exist $x^*\in\partial\,D_{K_y}$ and $\lambda^*\in(0,\bar{\lambda})$ such that $x^*=y+\lambda^*\, T(x^*)$. 
\end{theorem}

\noindent
{\bf Proof.}
If $T_{(y,\lambda)}$ has a fixed point on $\partial\,D_{K_y}$ for some $\lambda\in(0,\bar{\lambda})$ we are done. Otherwise, suppose that $T_{(y,\lambda)}$ is fixed point free on $\partial\,D_{K_y}$.
Now, observe that for each $\lambda\in(0,\bar{\lambda}]$ the operator $T_{(y,\lambda)}:\overline{D}_{K_y} \longrightarrow K_y$ satisfies condition \eqref{cond_Ty} and that $T_{(y,\lambda)}\left(\overline{D}_{K_y}\right)$ is relatively compact, which implies that the fixed point index $i_{K_y}(T_{(y,\lambda)},D_{K_y})$ is well--defined according to Definition \ref{def_index}.

Consider the map $H:\overline{D}_{K_y}\times[0,1]\rightarrow K_y$ defined as
\[H(x,t)=y+t\,\bar{\lambda}\,T(x). \]
Note that $H$ satisfies conditions $(a)$--$(c)$ in Theorem \ref{th_prop_index}, $i$. Hence, if $x\neq H(x,t)$ for all $(x,t)\in\partial\,D_{K_y}\times\left[0,1\right]$, then 
\[i_{K_y}(T_{(y,\bar{\lambda})},D_{K_y})=i_{K_y}(H(\cdot,1),D_{K_y})=i_{K_y}(H(\cdot,0),D_{K_y})=i_{K_y}(y,D_{K_y}). \]
By the normalization property, since $y\in D_{K_y}$, we have
\[i_{K_y}(T_{(y,\bar{\lambda})},D_{K_y})=i_{K_y}(y,D_{K_y})=1,\]
a contradiction. In conclusion, there exist $t^*\in(0,1)$ and $x^*\in\partial\,D_{K_y}$ such that $x^*=y+t^*\bar{\lambda}\, T(x^*)$.
\qed

\begin{remark} \label{rem1}
Note that Theorem \ref{th1_BK} has a two-fold interest: not only is a generalization of Theorem~\ref{th_BK_affine} to the context of discontinuous operators,
but also an improvement in the continuous case, since the conditions on the index are weakened and the result is extended to the setting of wedges.
\end{remark}

We now prove a result in the setting of normal cones which can be of a more direct applicability due to the use of the norm, as in the classical Birkhoff--Kellogg Theorem.

\begin{theorem}\label{th2_BK}
	Let $K\subset X$ be a normal cone with normal constant $c>0$ in a Banach space $X$, $D\subset X$ be an open bounded set with $0\in D$ and $y\in X$ be fixed. Assume that $T:\overline{D}_{K_y} \longrightarrow K$ is a mapping such that $T\,\overline{D}_{K_y}$ is relatively compact and
	\[\inf\{\left\|v\right\|:v\in\mathbb{T}x, \ x\in\partial\,D_{K_y} \}>0. \]
	
	If there exists a positive number 
	\[\bar{\lambda}>\dfrac{c\,\sup_{x\in\partial\,D}\left\|x\right\|}{\inf\{\left\|v\right\|:v\in\mathbb{T}x, \ x\in\partial\,D_{K_y} \}}, \]
	such that the operator $T_{(y,\lambda)}$ satisfies condition \eqref{cond_Ty} for each $\lambda\in(0,\bar{\lambda}]$,
	then there exist $x^*\in\partial\,D_{K_y}$ and $\lambda^*\in(0,\bar{\lambda})$ such that $x^*=y+\lambda^*\, T(x^*)$. 
\end{theorem}

\noindent
{\bf Proof.} We shall show that $i_{K_y}(T_{(y,\bar{\lambda})},D_{K_y})=0$ and so the conclusion is obtained as a consequence of Theorem \ref{th1_BK}.

Take $x_0\in K\setminus\{0\}$ and let us see that
\[x\notin y+\bar{\lambda}\,\mathbb{T}x+\beta\,x_0 \quad \text{ for all } \beta\geq 0 \text{ and } x\in\partial\,D_{K_y}. \]
Indeed, suppose that there exist $x_1\in\partial\,D_{K_y}$, $v\in\mathbb{T}x_1$ and $\beta_0\geq 0$ such that 
\[x_1=y+\bar{\lambda}\,v+\beta_0\, x_0. \]
Then $\bar{\lambda}\, v\preceq \bar{\lambda}\, v+\beta_0\, x_0=x_1-y$ and, since $K$ is normal,
\[\bar{\lambda}\left\|v\right\|\leq c\left\|x_1-y\right\|.\] 
Observe that $x_1-y\in\partial\, D$ and so
\[ \bar{\lambda}\left\|v\right\|\leq c\left(\sup_{x\in\partial\,D}\left\|x\right\| \right), \]
a contradiction with the choice of $\bar{\lambda}$.

On the other hand, since $\overline{D}_{K_y}$ and $\mathbb{T}\left(\overline{D}_{K_y} \right)$ are bounded, there exists $\bar{\beta}>0$ such that 
\[x\notin y+\bar{\lambda}\,\mathbb{T}x+\bar{\beta}\,x_0 \quad \text{ for all } x\in \overline{D}_{K_y}. \]
Consider the multivalued homotopy $\mathcal{H}:\overline{D}_{K_y}\times[0,1]\rightarrow 2^{K_y}$ defined as
\[\mathcal{H}(x,t)=y+\bar{\lambda}\,\mathbb{T}x+t\,\bar{\beta}\,x_0. \]
By the homotopy invariance property of the index for usc multivalued maps \cite{fp1},
\[i_{K_y}(\mathbb{T}_{(y,\bar{\lambda})},D_{K_y})=i_{K_y}(\mathcal{H}(\cdot,0),D_{K_y})=i_{K_y}(\mathcal{H}(\cdot,1),D_{K_y})=0. \]
Therefore, it follows from Proposition \ref{prop_indT=indTT} that $i_{K_y}(T_{(y,\bar{\lambda})},D_{K_y})=i_{K_y}(\mathbb{T}_{(y,\bar{\lambda})},D_{K_y})=0$.
\qed

The following corollary can be seen as an analogue of the classical result of
Krasnosel'ski\u{i} and Lady\v{z}enski\u{\i}.

\begin{corollary} \label{cor_BK}
	Let $K\subset X$ be a normal cone in a Banach space $X$ and $D\subset X$ be an open bounded set with $y\in D_{K_y}$. Assume that $T:\overline{D}_{K_y} \longrightarrow K$ is a mapping such that $T\,\overline{D}_{K_y}$ is relatively compact and, for each $\lambda>0$, the operator $T_{(y,\lambda)}$ satisfies condition \eqref{cond_Ty}.	
	If
	\[\inf\{\left\|v\right\|:v\in\mathbb{T}x, \ x\in\partial\,D_{K_y} \}>0, \]
	then there exist $x^*\in\partial\,D_{K_y}$ and $\lambda^*>0$ such that $x^*=y+\lambda^*\, T(x^*)$. 
\end{corollary}

\begin{remark}\label{rem2}
	Note that, in the non-affine case, Corollary \ref{cor_BK}	
	extends Theorem~\ref{th_KL} within the setting of discontinuous operators
	in normal cones. We stress that, in the non-affine case, Corollary \ref{cor_BK}
	can also be deduced as a consequence of the multivalued generalization of the Birkhoff--Kellogg theorem given in \cite{fp2}.
\end{remark}

\section{Applications} \label{sec-appl}

Consider the second order parameter-dependent differential equation
\begin{equation}\label{eq_2or}
	u''(t)+\lambda\,f(t,u(t),u(\sigma(t)))=0, \quad t\in[0,1],
\end{equation}
with initial conditions of the form
\begin{equation}\label{eq_initial}
	u(t)=\omega(t), \quad t\in[-r,0],
\end{equation}
and the final homogeneous boundary condition
\begin{equation}\label{BC}
	u(1)=0,
\end{equation}
where $\lambda$ is a positive parameter, $r\geq 0$, and $\sigma:[0,1]\rightarrow[-r,1]$ and $\omega:[-r,0]\rightarrow[0,\infty)$ are continuous functions. The nonlinearity $f:[0,1]\times[0,\infty)\times[0,\infty)\rightarrow[0,\infty)$ may be discontinuous with respect to the second argument in a sense which will be specified later.  

In order to study the problem \eqref{eq_2or}--\eqref{BC}, we shall use a superposition principle as in \cite{CI23}. To do so, first consider the Dirichlet BVP
\[\left\{\begin{array}{ll} u''(t)+y(t)=0, & \quad t\in[0,1], \\ u(0)=u(1)=0, \end{array} \right. \] 
whose unique solution is given by
\[u(t)=\int_{0}^{1}G(t,s)y(s)\,ds, \]
where $G$ is the corresponding Green's function. It is well-known that
\[G(t,s)=\left\{\begin{array}{ll} t(1-s), & \text{ if } 0\leq t\leq s\leq 1, \\ (1-t)s, & \text{ if } 0\leq s<t\leq 1, \end{array} \right. \]
and, moreover (see \cite{Lan-Webb}),
\begin{align*}
	G(t,s)&\leq \Phi(s), \quad t,s\in[0,1], \\ \dfrac{1}{4}\Phi(s)&\leq G(t,s), \quad t\in\left[\dfrac{1}{4},\dfrac{3}{4}\right], \ s\in[0,1],
\end{align*}
with $\Phi(s):=s(1-s)$. Associated to the Green's function, we consider the kernel $k:[-r,1]\times[0,1]\rightarrow\mathbb{R}$ defined as
\begin{equation}
k(t,s)=	\left\{\begin{array}{ll} G(t,s), & \text{ if } t\geq 0, \\ 0, & \text{ if } t<0. \end{array} \right.
\end{equation}

On the other hand, note that the function $\hat{y}(t)=1-t$ solves the Dirichlet BVP
\[\left\{\begin{array}{ll} u''(t)=0, & \quad t\in[0,1], \\ u(0)=1, \ u(1)=0, \end{array} \right. \]
so we define the function
\begin{equation}\label{eq_y}
y(t)=\left\{\begin{array}{ll} \omega(t), & \text{ if } t\leq 0, \\ \hat{y}(t)\,\omega(0), & \text{ if } t>0, \end{array} \right. 
\end{equation}
which will be the \textit{vertex} of our affine cone.

In order to apply the theory of the previous section, we will work in the Banach space of continuous functions $X=\mathcal{C}([-r,1])$, endowed with the usual sup-norm, $\left\|\cdot\right\|_{[-r,1]}$, and the cone
\[K=\left\{u\in\mathcal{C}([-r,1],[0,\infty)):u(t)=0 \text{ for all } t\in[-r,0], \min_{t\in[1/4,3/4]}u(t)\geq \dfrac{1}{4}\left\|u\right\|_{[0,1]} \right\}. \]
Observe that $K$ is a normal cone with normal constant $c=1$ and that $\left\|u\right\|_{[-r,1]}=\left\|u\right\|_{[0,1]}$ for all $u\in K$. Now, for the vertex $y$ defined in \eqref{eq_y}, we consider the translate of the cone $K$ given by
\[K_y:=y+K=\{y+u:u\in K \}, \]
and, for each $\rho>0$, we denote by $K_{y,\rho}$ the relatively open bounded set 
\[K_{y,\rho}:=\left\{y+u:u\in K, \ \left\|u\right\|_{[0,1]}<\rho \right\}. \]

We will look for solutions of the following perturbed Hammerstein integral equation
\begin{equation}\label{eq_H}
	u(t)=y(t)+\lambda\int_{0}^{1}k(t,s)f(s,u(s),u(\sigma(s)))\,ds=:y(t)+\lambda\,Tu(t), \quad t\in[-r,1],
\end{equation}
located in the affine cone $K_y$. 

\begin{definition}
By a solution of  the problem \eqref{eq_2or}--\eqref{BC}
we mean a solution $u \in  \mathcal{C}([-r, 1], \R)$
of the integral equation \eqref{eq_H}.
\end{definition}

Before doing so, we need to define the type of regions where $f$ is allowed to be discontinuous. The concept of \textit{admissible discontinuity curve} used here has been widely employed in \cite{FLR}.

\begin{definition}
	A $\lambda$-admissible discontinuity curve for the second-order parameter dependent differential equation $u''+\lambda\,f(t,u,u(\sigma))=0$ is a $W^{2,1}$ function $\gamma:[a,b]\subset [0,1]\rightarrow [0,\infty)$ satisfying that there exist $\varepsilon>0$ and $\psi\in L^1(a,b)$, $\psi(t)>0$ for a.a. $t\in[a,b]$ such that either
	\begin{equation}\label{eq_ad1}
	-\gamma''(t)+\psi(t)<\lambda\,f(t,y,z) \quad \text{ for a.a. } t\in[a,b], \text{ all } y\in[\gamma(t)-\varepsilon,\gamma(t)+\varepsilon] \text{ and all } z\in[0,\infty) 
	\end{equation}
	or 
	\begin{equation}\label{eq_ad2}
	-\gamma''(t)-\psi(t)>\lambda\,f(t,y,z) \quad \text{ for a.a. } t\in[a,b], \text{ all } y\in[\gamma(t)-\varepsilon,\gamma(t)+\varepsilon] \text{ and all } z\in[0,\infty). 
	\end{equation}	
\end{definition}

\begin{remark}\label{rmk_gamma}
	Since $f$ is non-negative, in order to have that $\gamma$ is a $\lambda$-admissible discontinuity curve for the differential equation \eqref{eq_2or} and any $\lambda>0$, it suffices that 
	\[0<\gamma''(t) \quad \text{ for a.a. } t\in[a,b]. \] 
	Indeed, one may check that condition \eqref{eq_ad1} holds with $\psi(t)=\gamma''(t)/2$, $t\in[a,b]$.
\end{remark}

Let us now state and prove the main result of this Section.

\begin{theorem}\label{th}
	Let $\rho>0$ and assume that the following conditions hold:
	\begin{enumerate}
		\item[$(H_1)$] any composition $t\in [0,1]\mapsto f(t,u(t),v(t))$ is measurable provided that $u,v\in \mathcal{C}([0,1],[0,\infty))$;
		\item[$(H_2)$] there exists $M_{\rho}\in L^1([0,1])$ such that
		\[f(t,u,v)\leq M_{\rho}(t) \quad \text{ for a.a. } t\in[0,1] \text{ and all } (u,v) \text{ with } 0\leq u,v \leq \rho+\left\|\omega\right\|_{[-r,0]}; \]
		\item[$(H_3)$] there exists $\delta_{\rho}\in L^1([1/4,3/4])$ such that
		\[f(t,u,v)\geq \delta_{\rho}(t) \quad \text{ for a.a. } t\in\left[1/4,3/4\right] \text{ and all } (u,v) \text{ with } 0\leq u,v \leq \rho+\left\|\omega\right\|_{[-r,0]} \]
		and 
		\[\bar{\delta}:=\sup_{t\in[1/4,3/4]}\int_{1/4}^{3/4}k(t,s)\delta_{\rho}(s)\,ds>0; \]
		\item[$(H_4)$] there exists a countable number of curves $\gamma_n:I_n=[a_n,b_n]\rightarrow [0,\infty)$, $n\in\mathbb{N}$, such that for a.a. $t\in [0,1]$ the function $f(t,\cdot,\cdot)$ is continuous on $\left([0,\infty)\setminus \bigcup_{n:t\in I_n}\{\gamma_n(t) \} \right)\times [0,\infty)$ and, moreover, each $\gamma_n$ is a $\lambda$-admissible discontinuity curve for each $\lambda\in(0,\bar{\lambda}]$ and some $\bar{\lambda}>\rho/\bar{\delta}$.
	\end{enumerate}
	Then there exist $\lambda_{\rho}\in(0,\bar{\lambda})$ and $u_{\rho}\in \partial \,K_{y,\rho}$ that satisfy the integral equation \eqref{eq_H}.
\end{theorem}

\noindent
{\bf Proof.} Let us divide the proof in several steps:

\textbf{Step 1.} \textit{The operator $T$, defined in \eqref{eq_H}, maps the set $\overline{K}_{y,\rho}$ into the cone $K$ and, moreover, $T\,\overline{K}_{y,\rho}$ is relatively compact.} 

First, let $u\in \overline{K}_{y,\rho}$ be arbitrarily fixed and let us show that $T u\in K$. By definition,
\[Tu(t)=\int_{0}^{1}k(t,s)f(s,u(s),u(\sigma(s)))\,ds, \quad t\in[-r,1]. \]
 The continuity of the kernel $k$, jointly with hypothesis $(H_1)$, $(H_2)$ and the constant sign of $f$ and $k$, imply that $T u\in \mathcal{C}([-r,1],[0,\infty))$. Moreover, since $k(t,s)=0$ for all $t\leq 0$, we have that $Tu(t)=0$ for all $t\in[-r,0]$. Now, for $t\in[1/4,3/4]$, we have
 \begin{align*}
 	Tu(t)&=\int_{0}^{1}G(t,s)f(s,u(s),u(\sigma(s)))\,ds \\ &\geq \dfrac{1}{4}\int_{0}^{1}\Phi(s)f(s,u(s),u(\sigma(s)))\,ds \geq \dfrac{1}{4}\left\|Tu \right\|_{[0,1]},
 \end{align*}
 as a consequence of the properties of the Green's function $G$ stated above.
 In conclusion, $Tu\in K$.

On the other hand, the compactness of the set $\overline{T\,\overline{K}_{y,\rho}}$ follows from assumption $(H_2)$ and the continuity of the kernel $k$, combined with a careful use of the Arzel\`{a}--Ascoli theorem (see \cite{Webb-Cpt}).

\textbf{Step 2.} \textit{For each $\lambda\in(0,\bar{\lambda}]$, the operator $y+\lambda\,T$ satisfies that
	\begin{equation}\label{cond_TT_proof}
		\left\{u\right\}\cap\left\{y+\lambda\,\mathbb{T}(u) \right\}\subset\left\{y+\lambda\,T(u)\right\} \quad \mbox{for every $u\in\overline{K}_{y,\rho}$,}	
	\end{equation}
	where $\bar{\lambda}$ is fixed by hypothesis $(H_4)$.}

Fix arbitrary $\lambda\in(0,\bar{\lambda}]$ and $u\in\overline{K}_{y,\rho}$. Now, consider two different cases:

\textit{Case 1.} $m\left(\left\{t\in I_n:u(t)=\gamma_n(t) \right\}\right)=0$ for all $n\in\mathbb{N}$ (where $m$ denotes Lebesgue measure). 

Let us prove that $T$ is continuous at $u$, which implies that $\mathbb{T}(u)=\left\{T(u) \right\}$ and thus condition \eqref{cond_TT_proof} holds for such $u$. Indeed, in this case we have that for a.a. $t\in[0,1]$ the function $f(t,\cdot,\cdot)$ is continuous at $(u(t),u(\sigma(t)))$. Hence, if $u_k\to u$ uniformly in $[-r,1]$, then
\[f(t,u_k(t),u_k(\sigma(t)))\to f(t,u(t),u(\sigma(t))) \quad \text{for a.a. } t\in[0,1], \] 
which implies, due to Lebesgue's dominated convergence theorem, that $T u_k\to Tu$ in $\mathcal{C}([-r,1])$.

\textit{Case 2.} $m\left(\left\{t\in I_n:u(t)=\gamma_n(t) \right\}\right)>0$ for some $n\in\mathbb{N}$.

In this case, one can show that $u\notin y+\lambda\,\mathbb{T}(u)$, which implies that condition \eqref{cond_TT_proof} holds for such $u$. The proof is based on condition $(H_4)$ and the fact that the function $\gamma_n$ is a $\lambda$-admissible discontinuity curve for the problem. It can be replicated following the reasoning in the proof of Proposition 4.7, Case 2, in \cite{FLR_JFPTA}.

\textbf{Step 3.} \textit{It holds that 
	\[\inf\left\{\left\|v\right\|_{[-r,1]}:v\in\mathbb{T}u, \ u\in\partial\,K_{y,\rho} \right\}\geq \bar{\delta}>0.\]}

For $u\in\partial\,K_{y,\rho}$ and $\varepsilon>0$, take $u_i\in \overline{B}_{\varepsilon}(u)\cap \overline{K}_{y,\rho}$ and $\lambda_i\geq 0$ with $\sum \lambda_i=1$, $i=1,2,\dots,m$. Then, by assumption $(H_3)$, we have for $t\in[1/4,3/4]$,
\begin{align*}
	\sum_{i=1}^{m}\lambda_i\,Tu_i(t)&=\sum_{i=1}^{m}\lambda_i\,\int_{0}^{1}k(t,s)f(s,u_i(s),u_i(\sigma(s)))\,ds \\ &\geq \sum_{i=1}^{m}\lambda_i\,\int_{1/4}^{3/4}k(t,s)f(s,u_i(s),u_i(\sigma(s)))\,ds \\ &\geq \sum_{i=1}^{m}\lambda_i\,\int_{1/4}^{3/4}k(t,s)\delta_{\rho}(s)\,ds \\ &=\int_{1/4}^{3/4}k(t,s)\delta_{\rho}(s)\,ds.
\end{align*}
Hence, for any $v\in {\rm co} \, T\left(\overline{B}_{\varepsilon}(u)\cap \overline{K}_{y,\rho}\right)$, we have
\[\left\|v\right\|_{[-r,1]}\geq \left\|v\right\|_{[1/4,3/4]}\geq \sup_{t\in[1/4,3/4]}\int_{1/4}^{3/4}k(t,s)\delta_{\rho}(s)\,ds=\bar{\delta}. \]
Since $\mathbb{T} u\subset \overline{\rm co} \, T\left(\overline{B}_{\varepsilon}(u)\cap \overline{K}_{y,\rho}\right)$, it follows that 
$\left\|v\right\|_{[-r,1]}\geq \bar{\delta}>0$ for any $v\in \mathbb{T} u$, as wished.

\medskip

Therefore, the conclusion follows from Theorem \ref{th2_BK}.
\qed

\begin{remark}
We emphasize that hypotheses $(H_1)$, $(H_2)$ and $(H_4)$ do not imply that $f$ be a Carath\'eodory map, since due to $(H_4)$ the function $f$ can be discontinuous with respect to the last variables.
	
	Furthermore, note that if, for each $(x,y)\in[0,\infty)\times [0,\infty)$, the map $t\in[0,1]\mapsto f(t,x,y)$ is measurable and, for a.a. $t\in[0,1]$, the map $(x,y)\mapsto f(t,x,y)$ is continuous, then condition $(H_1)$ holds. However, the measurability of the map  $t\in[0,1]\mapsto f(t,x,y)$ together with $(H_4)$ do not imply necessarily that condition $(H_1)$ holds. More information about the measurability of compositions in this setting can be found in \cite[Section 3.1]{CidPo}.
\end{remark}

\begin{corollary}
	Let $\rho>0$ and assume that conditions $(H_1)$--$(H_3)$ hold and, moreover,
	\begin{enumerate}
		\item[$(H_4^*)$] there exists a countable number of curves $\gamma_n:I_n=[a_n,b_n]\rightarrow [0,\infty)$, $n\in\mathbb{N}$, such that $\gamma_n''>0$ and for a.a. $t\in [0,1]$ the function $f(t,\cdot,\cdot)$ is continuous on $\left([0,\infty)\setminus \bigcup_{n:t\in I_n}\{\gamma_n(t) \} \right)\times [0,\infty)$.
	\end{enumerate}
Then there exist $\lambda_{\rho}>0$ and $u_{\rho}\in \partial \,K_{y,\rho}$ that satisfy the integral equation \eqref{eq_H}.
\end{corollary}

\noindent
{\bf Proof.} It follows from Theorem \ref{th} together with Remark \ref{rmk_gamma}.
\qed

Consider the special case of \eqref{eq_2or} where the nonlinearity can be seen as a discontinuous perturbation of a Carath\'eodory function, that is,
\begin{equation}\label{eq_2or_g+h}
	u''(t)+\lambda\left(g(t,u(\sigma(t)))+h(u(t)) \right)=0, \quad t\in[0,1],
\end{equation}
where $g:[0,1]\times[0,\infty)\to [0,\infty)$ is a Carath\'eodory function and $h:[0,\infty)\to[0,\infty)$ is locally bounded and continuous except at most at a countable number of points.

\begin{corollary}
	Assume that the following conditions hold:
	\begin{enumerate}
		\item[$(C_1)$] $g$ satisfies the Carath\'eodory conditions, namely,
		\begin{enumerate}
			\item $g(\cdot,v)$ is measurable for each fixed $v\in[0,\infty)$;
			\item $g(t,\cdot)$ is continuous for a.a. $t\in[0,1]$;
			\item for each $R>0$, there exists $M_R\in L^1([0,1])$ such that
			\[g(t,v)\leq M_R(t) \quad \text{ for a.a. } t\in[0,1] \text{ and all } v\in[0,R]; \]
		\end{enumerate} 
		\item[$(C_2)$] $h$ is locally bounded, $t\mapsto h(u(t))$ is measurable for each non-negative continuous function $u$
		and there exists a countable set $A$ such that $h$ is continuous in $[0,\infty)\setminus A$;
		\item[$(C_3)$] there exists $\delta\in L^1([0,1])$, $\delta(t)>0$ for a.a. $t\in[0,1]$ such that 
		\[g(t,v)\geq \delta(t) \quad \text{ for a.a. } t\in[0,1] \text{ and all } v\geq 0. \]
	\end{enumerate}
	Then for each $\rho>0$ there exist $\lambda_{\rho}>0$ and $u_{\rho}\in \partial \,K_{y,\rho}$ that satisfy the BVP \eqref{eq_2or_g+h}--\eqref{eq_initial}--\eqref{BC}.
\end{corollary}

\noindent
{\bf Proof.} Observe that, for each $\rho>0$, Theorem \ref{th} can be applied to the function
\[f(t,u,v)=g(t,v)+h(u). \]
Note that hypotheses $(H_1)$--$(H_3)$ are satisfied. 

Now, consider the countable set $A$ where $h$ may be discontinuous and denote $A=\{a_k:k\in\mathbb{N} \}$. Define the constant functions $\gamma_n:[0,1]\to [0,\infty)$ given by $\gamma_n(t)=a_n$, $t\in[0,1]$, $n\in\mathbb{N}$. For each $\lambda>0$ fixed, choose the $L^1$--function $\psi(t)=\lambda\,\delta(t)$, $t\in[0,1]$. Then each function $\gamma_n$ satisfies condition \eqref{eq_ad1} and so it is a $\lambda$-admissible discontinuity curve. 
\qed

Now, let us restrict our efforts to the particular case of problem \eqref{eq_2or}--\eqref{BC} in which the deviated argument is given by a continuously differentiable function $\sigma$ with constant derivative equal to $1$ or $-1$. Notice that it covers the meaningful situations of equations with delay (where $\sigma(t)=t-r$) or with reflection of the argument (where, for instance, $\sigma(t)=1-r-t$). 

In this case, we are able to prove another version of Theorem \ref{th} where the nonlinearity $f$ may be \textit{more} discontinuous. More precisely, we weaken assumption $(H_4)$ allowing $f$ to be discontinuous w.r.t. the second and third variables over the graphs of two countable families of functions.

\begin{theorem}\label{th_D}
	Let $\rho>0$ and assume conditions $(H_1)$--$(H_3)$ in Theorem \ref{th} hold. Moreover, suppose that $\sigma:[0,1]\rightarrow[-r,1]$ is a continuously differentiable function with constant derivative $\sigma'=\pm 1$ and the following assumption holds:
	\begin{itemize}
		\item[($D$)] there exist two countable families of curves $\gamma_n:I_n=[a_n,b_n]\rightarrow [0,\infty)$, $n\in\mathbb{N}$, and $\Gamma_j:\mathcal{I}_j=[c_j,d_j]\rightarrow [0,\infty)$, $j\in\mathbb{N}$, such that for a.a. $t\in [0,1]$ the function 
		\[f(t,\cdot,\cdot) \text{ is continuous on } \left([0,\infty)\setminus \bigcup_{n:t\in I_n}\{\gamma_n(t) \} \right)\times \left([0,\infty)\setminus \bigcup_{j:t\in \mathcal{I}_j}\{\Gamma_j(t) \} \right).\] 
		For each $\lambda\in(0,\bar{\lambda}]$ with $\bar{\lambda}>\rho/\bar{\delta}$, each function $\gamma_n$ is a $\lambda$-admissible discontinuity curve and each function $\Gamma_j$ satisfies that
		\begin{enumerate}
			\item[$(a)$] $\Gamma_j(t)\neq \omega(\sigma(t))$ for a.a. $t\in \mathcal{I}_j\cap \sigma^{-1}([-r,0])$;  
			\item[$(b)$] the restriction of $\Gamma_j$ to $\mathcal{I}_j\cap \sigma^{-1}([0,1])$ satisfies either of the following conditions: there exists $\varepsilon_j>0$ and $\psi_j\in L^1(\mathcal{I}_j)$, $\psi_j(t)>0$ for a.a. $t\in \mathcal{I}_j\cap \sigma^{-1}([0,1])$ such that
			\begin{enumerate}
				\item $-\Gamma_j''(t)+\psi_j(t)<\lambda\,f(\sigma(t),y,z)$ for a.a. $t\in \mathcal{I}_j\cap\sigma^{-1}([0,1])$,  all  $y\in[\Gamma_j(t)-\varepsilon_j,\Gamma_j(t)+\varepsilon_j]$  and all $z\in[0,\infty)$; or
				\item $-\Gamma_j''(t)-\psi_j(t)>\lambda\,f(\sigma(t),y,z)$ for a.a. $t\in \mathcal{I}_j\cap\sigma^{-1}([0,1])$,  all  $y\in[\Gamma_j(t)-\varepsilon_j,\Gamma_j(t)+\varepsilon_j]$  and all $z\in[0,\infty)$
			\end{enumerate}
		\end{enumerate}
	\end{itemize} 

	Then there exist $\lambda_{\rho}\in(0,\bar{\lambda})$ and $u_{\rho}\in \partial \,K_{y,\rho}$ that satisfy the integral equation \eqref{eq_H}, that is, they solve the problem \eqref{eq_2or}--\eqref{BC}. 
\end{theorem}

\noindent
{\bf Proof.} It follows in line of the proof of Theorem \ref{th} as a consequence of Theorem \ref{th2_BK}. Observe that it suffices to rewrite Step 2. Let us prove that for each $\lambda\in(0,\bar{\lambda}]$, the operator $y+\lambda\,T$ satisfies that
\begin{equation}\label{cond_TT_proof2}
	\left\{u\right\}\cap\left\{y+\lambda\,\mathbb{T}(u) \right\}\subset\left\{y+\lambda\,T(u)\right\} \quad \mbox{for every $u\in\overline{K}_{y,\rho}$,}	
\end{equation}
where $\bar{\lambda}$ is fixed by assumption $(D)$.

Fix arbitrary $\lambda\in(0,\bar{\lambda}]$ and $u\in\overline{K}_{y,\rho}$. Now, consider three different cases:

\textit{Case 1.} $m\left(\left\{t\in I_n:u(t)=\gamma_n(t) \right\}\cup \left\{t\in \mathcal{I}_j:u(\sigma(t))=\Gamma_j(t) \right\}\right)=0$ for all $j,n\in\mathbb{N}$. Then for a.a. $t\in[0,1]$ the function $f(t,\cdot,\cdot)$ is continuous at $(u(t),u(\sigma(t)))$ and so $T$ is continuous at $u$.

\textit{Case 2.} $m\left(\left\{t\in \mathcal{I}_j:u(\sigma(t))=\Gamma_j(t) \right\}\right)>0$ for some $j\in\mathbb{N}$. Let us prove that $u\notin y +\lambda\,\mathbb{T} u$, which can be justified as in the proof of \cite[Proposition 4.7]{FLR_JFPTA}, but we include the reasoning here again for completeness.

Since $u\in K_y$, we have $u(s)=y(s)=\omega(s)$ for all $s\in[-r,0]$ and thus $u(\sigma(t))=\omega(\sigma(t))$ for all $t\in \sigma^{-1}([-r,0])$. Now condition $(D)$, $(a)$, implies that
\[m\left(\left\{t\in \mathcal{I}_j\cap \sigma^{-1}([-r,0]):u(\sigma(t))=\Gamma_j(t) \right\}\right)=0 \quad \text{for all } j\in\mathbb{N}. \]
Hence, we can fix some $j\in\mathbb{N}$ such that
\[m\left(\left\{t\in \mathcal{I}_j\cap \sigma^{-1}([0,1]):u(\sigma(t))=\Gamma_j(t) \right\}\right)>0. \]
By condition $(D)$, $(b)$, we can assume that there exist $\varepsilon_j>0$ and $\psi_j\in L^1(\mathcal{I}_j)$, $\psi_j(t)>0$ on $\mathcal{I}_j$, such that
\begin{equation}\label{eq_Gamma}
-\Gamma_j''(t)+\psi_j(t)<\lambda\,f(\sigma(t),y,z) \quad \text{for a.a. } t\in \mathcal{I}_j\cap\sigma^{-1}([0,1]), \text{ all } y\in[\Gamma_j(t)-\varepsilon_j,\Gamma_j(t)+\varepsilon_j] \text{ and all } z\in[0,\infty). 
\end{equation}
In what follows, let us denote $J:=\left\{t\in \mathcal{I}_j\cap \sigma^{-1}([0,1]):u(\sigma(t))=\Gamma_j(t) \right\}$ and $M:=\lambda\left(M_{\rho}\circ\sigma\right)$. By technical results of Lebesgue measure (see \cite[Lemma 4.2 and Corollary 4.3]{FLR_JFPTA}), we know that there exists a measurable set $J_0\subset J$ with $m(J)=m(J_0)$ such that for all $\tau_0\in J_0$, 
\begin{equation}\label{eq_J0}
	\lim_{t\to \tau_0^+}\dfrac{\int_{[\tau_0,t]\setminus J}M(s)\,ds}{\int_{\tau_0}^{t}\psi_j(s)\,ds}=0=\lim_{t\to \tau_0^-}\dfrac{\int_{[t,\tau_0]\setminus J}M(s)\,ds}{\int_{t}^{\tau_0}\psi_j(s)\,ds}
\end{equation}
and, moreover, there is $J_1\subset J_0$ with $m(J_0\setminus J_1)=0$ such that for all $\tau_0\in J_1$,
\begin{equation}\label{eq_J1}
	\lim_{t\to \tau_0^+}\dfrac{\int_{[\tau_0,t]\cap J_0}\psi_j(s)\,ds}{\int_{\tau_0}^{t}\psi_j(s)\,ds}=1=\lim_{t\to \tau_0^-}\dfrac{\int_{[t,\tau_0]\cap J_0}\psi_j(s)\,ds}{\int_{t}^{\tau_0}\psi_j(s)\,ds}.
\end{equation}
Now, fix $\tau_0\in J_1$. By \eqref{eq_J0} and \eqref{eq_J1}, there exists $\bar{t}>0$ sufficiently close to $0$ such that for all $t\in[\bar{t},2\,\bar{t}]$ the following inequalities hold
\begin{align*}
	2\int_{[\tau_0,\tau_0+t]\setminus J}M(s)\,ds<\dfrac{1}{4}\int_{\tau_0}^{\tau_0+t}\psi_j(s)\,ds, \\ 
	\int_{[\tau_0,\tau_0+t]\cap J}\psi_j(s)\,ds>\dfrac{1}{2}\int_{\tau_0}^{\tau_0+t}\psi_j(s)\,ds \\ 
	\intertext{and} 
	2\int_{[\tau_0-t,\tau_0]\setminus J}M(s)\,ds<\dfrac{1}{4}\int_{\tau_0-t}^{\tau_0}\psi_j(s)\,ds, \\ 
	\int_{[\tau_0-t,\tau_0]\cap J}\psi_j(s)\,ds>\dfrac{1}{2}\int_{\tau_0-t}^{\tau_0}\psi_j(s)\,ds.
\end{align*}
Define the positive number
\[r=\dfrac{\bar{t}}{8}\min\left\{\int_{\tau_0-\bar{t}}^{\tau_0}\psi_j(s)\,ds,\int_{\tau_0}^{\tau_0+\bar{t}}\psi_j(s)\,ds \right\}. \]

Let $\varepsilon_j>0$ be given above. Let us show that for every finite family $u_i\in\overline{B}_{\varepsilon_j}(u)\cap \overline{K}_{y,\rho}$ and $\mu_i\in[0,1]$ ($i=1,2,\dots,m$), with $\sum\mu_i=1$, we have
\begin{equation}\label{eq_r} 
\left\|u-\left(y+\lambda\sum_{i=1}^{m}\mu_i Tu_i\right)\right\|_{[0,1]}\geq r, 
\end{equation}
which implies $u\notin y +\lambda\,\mathbb{T} u$. Notice that we can suppose without loss of generality that the restriction of $u$ to $[0,1]$, denoted also as $u$, satisfies that $u\in y+\lambda\,Q$ where $Q$ is the subset of $\mathcal{C}([0,1])$ defined as
\[Q=\left\{u\in\mathcal{C}^1([0,1]):\left|u'(t)-u'(s)\right|\leq \int_{s}^{t}M_{\rho}(r)\,dr \ \text{ whenever } 0\leq s\leq t\leq 1 \right\}. \]
Indeed, due to assumption $(H_2)$, then 
\[T u(t)=\int_{0}^{1}k(t,s)f(s,u(s),u(\sigma(s)))\,ds\leq \int_{0}^{1}M_{\rho}(s)\,ds. \]
Therefore, $T\left(\overline{K}_{y,\rho}\right)\subset Q$ and $Q$ is a closed convex subset of $\mathcal{C}([0,1])$ (see \cite[Lemma 4.5]{FLR_JFPTA}), which implies that $\mathbb{T}\left(\overline{K}_{y,\rho}\right)\subset Q$.

In order to prove \eqref{eq_r}, for simplicity, let us denote $z=\lambda\sum_{i=1}^{m}\mu_i Tu_i$ and $v=u-y$. For a.a. $t\in J$, we have by the chain rule that
\[\left(z\circ\sigma\right)''(t)=z''(\sigma(t))\left(\sigma'(t)\right)^2+z'(\sigma(t))\sigma''(t) \]
and, since $\sigma'=\pm 1$,
\[\left(z\circ\sigma\right)''(t)=z''(\sigma(t))=-\lambda\sum_{i=1}^{m}\mu_i f(\sigma(t),u_i(\sigma(t)),u_i(\sigma(\sigma(t)))). \]

On the other hand, for every $i\in\{1,\dots,m \}$ and $t\in J$, we deduce from $u_i\in \overline{B}_{\varepsilon_j}(u)$ that
\[\left|u_i(\sigma(t))-\Gamma_j(t) \right|=\left|u_i(\sigma(t))-u(\sigma(t))\right|\leq \varepsilon_j \]
and then condition \eqref{eq_Gamma} ensures that for a.a. $t\in J$,
\begin{align*}
	z''(\sigma(t))&=-\lambda\sum_{i=1}^{m}\mu_i f(\sigma(t),u_i(\sigma(t)),u_i(\sigma(\sigma(t)))) \\ &<\sum_{i=1}^{m}\mu_i\left(\Gamma_j''(t)-\psi_j(t)\right)=\Gamma_j''(t)-\psi_j(t)=u''(\sigma(t))-\psi_j(t).
\end{align*}
Note that $y''(s)=0$ for all $s>0$, so we obtain that for a.a. $t\in J$,
\[z''(\sigma(t))<u''(\sigma(t))-y''(\sigma(t))-\psi_j(t)=v''(\sigma(t))-\psi_j(t). \]
By integration, for $t\in[\bar{t},2\,\bar{t}]$,
\begin{align*}
	z'(\sigma(\tau_0))-z'(\sigma(\tau_0-t))&=\int_{\tau_0-t}^{\tau_0}z''(\sigma(s))\,ds=\int_{[\tau_0-t,\tau_0]\cap J}z''(\sigma(s))\,ds+\int_{[\tau_0-t,\tau_0]\setminus J}z''(\sigma(s))\,ds \\ &<\int_{[\tau_0-t,\tau_0]\cap J}v''(\sigma(s))\,ds-\int_{[\tau_0-t,\tau_0]\cap J}\psi_j(s)\,ds +\int_{[\tau_0-t,\tau_0]\setminus J}M(s)\,ds \\
	&=v'(\sigma(\tau_0))-v'(\sigma(\tau_0-t))-\int_{[\tau_0-t,\tau_0]\setminus J}v''(\sigma(s))\,ds -\int_{[\tau_0-t,\tau_0]\cap J}\psi_j(s)\,ds \\ &\quad +\int_{[\tau_0-t,\tau_0]\setminus J}M(s)\,ds \\ &\leq v'(\sigma(\tau_0))-v'(\sigma(\tau_0-t))-\int_{[\tau_0-t,\tau_0]\cap J}\psi_j(s)\,ds+2\int_{[\tau_0-t,\tau_0]\setminus J}M(s)\,ds \\ &<v'(\sigma(\tau_0))-v'(\sigma(\tau_0-t))-\dfrac{1}{4}\int_{\tau_0-t}^{\tau_0}\psi_j(s)\,ds.
\end{align*}
Hence, we have for all $t\in [\bar{t},2\,\bar{t}]$ that
\[z'(\sigma(\tau_0-t))-v'(\sigma(\tau_0-t))>z'(\sigma(\tau_0))-v'(\sigma(\tau_0))+\dfrac{1}{4}\int_{\tau_0-t}^{\tau_0}\psi_j(s)\,ds. \]

In case $z'(\sigma(\tau_0))\geq v'(\sigma(\tau_0))$, then
\[z'(\sigma(\tau_0-t))-v'(\sigma(\tau_0-t))>\dfrac{1}{4}\int_{\tau_0-\bar{t}}^{\tau_0}\psi_j(s)\,ds \quad \text{for all } t\in[\bar{t},2\,\bar{t}], \]
so, by integration,
\begin{align*} 
z(\sigma(\tau_0-\bar{t}))-v(\sigma(\tau_0-\bar{t}))&>z(\sigma(\tau_0-2\,\bar{t}))-v(\sigma(\tau_0-2\,\bar{t}))+\dfrac{\bar{t}}{4}\int_{\tau_0-\bar{t}}^{\tau_0}\psi_j(s)\,ds \\ &\geq z(\sigma(\tau_0-2\,\bar{t}))-v(\sigma(\tau_0-2\,\bar{t}))+2\,r.
\end{align*}
Then, either $z(\sigma(\tau_0-2\,\bar{t}))-v(\sigma(\tau_0-2\,\bar{t}))<-r$ or $z(\sigma(\tau_0-\bar{t}))-v(\sigma(\tau_0-\bar{t}))>r$, and thus $\left\|v-z\right\|>r$, as wished.

It can be seen in a similar way that
\[z'(\sigma(\tau_0))-v'(\sigma(\tau_0))>z'(\sigma(\tau_0+t))-v'(\sigma(\tau_0+t))+\dfrac{1}{4}\int_{\tau_0}^{\tau_0+t}\psi_j(s)\,ds, \quad t\in[\bar{t},2\,\bar{t}], \]
which ensures that $\left\|v-z\right\|>r$ if $z'(\sigma(\tau_0))< v'(\sigma(\tau_0))$.

\textit{Case 3.} $m\left(\left\{t\in I_n:u(t)=\gamma_n(t) \right\}\right)>0$ for some $n\in\mathbb{N}$. It follows as in Case 2.

Finally, Theorem \ref{th2_BK} gives the conclusion.
\qed

\begin{example} \label{example2}
	Consider the function $\phi:\mathbb{R}\rightarrow\mathbb{R}$ given by
	\[\phi(x)=\sum_{n:q_n<x} 2^{-n}, \]
	where $\{q_n \}_{n\in\mathbb{N}}$ is an enumeration of the rational numbers. Observe that $\phi$ is discontinuous at the rational numbers and continuous at the irrational ones. 
	
	We study the existence of solutions for the following BVP with delay
	\begin{equation}\label{eq_ex} 
	\left\{\begin{array}{rl} -u''(t)&=\lambda\left(\dfrac{2-\phi(u(t)-t^2)}{\sqrt{t}}+\phi(u(t-1/2)-t^2)(u(t))^3 \right), \quad t\in[0,1], \\[0.25cm] u(t)&=\sqrt{1+2\,t}, \quad t\in[-1/2,0], \\ u(1)&=0. \end{array} \right. 
	\end{equation}
	In this case,
	\[f(t,u,v)=\dfrac{2-\phi(u-t^2)}{\sqrt{t}}+\phi(v-t^2)u^3, \quad (t,u,v)\in[0,1]\times[0,+\infty)\times [0,+\infty).\]

	For a fixed $\rho>0$, we can choose 
	\[M_{\rho}(t)=\dfrac{2}{\sqrt{t}}+(\rho+1)^3 \quad \text{and} \quad \delta_{\rho}(t)=\delta(t)=\dfrac{1}{\sqrt{t}} \quad (t\in[0,1]), \]
	in order to check hypotheses $(H_2)$ and $(H_3)$. On the other hand, for each rational $q_n$ we define the function $\gamma_n:[0,1]\to\mathbb{R}$ as
	\[\gamma_n(t)=t^2+q_n, \]
	and so for a.a. $t\in[0,1]$, $f(t,\cdot,\cdot)$ is continuous on $\left([0,\infty)\setminus\bigcup_n \{\gamma_n(t) \} \right)\times \left([0,\infty)\setminus\bigcup_n \{\gamma_n(t) \} \right)$. Note that for each $n\in\mathbb{N}$ and each $\lambda>0$,
	\[-\gamma_n''(t)=-2<0\leq \lambda\, f(t,u,v) \quad \text{ for a.a. } t\in[0,1] \text{ and all } u,v\in[0,+\infty) \]
	and thus condition $(D)$ in Theorem \ref{th_D} holds. Therefore, this result ensures that the BVP \eqref{eq_ex} has uncountable many pairs of solutions and parameters $(u_{\rho},\lambda_{\rho})$.
\end{example}
\begin{remark}\label{eigen-clas}
We stress that the theory presented so far is applicable and represents a novelty even in the \emph{special case} of eigenvalue problems for ODEs, in presence of discontinuities. To illustrate this fact one may consider the eigenvalue problem 
\begin{equation}\label{eigendir}
u''(t)=\lambda \tilde{f}(u(t)),\, t \in (0,1);\ u(0)=u(1)=0,
\end{equation}
a classical problem studied in the book of Guo and Lakshmikantham~\cite[Example 2.3.2]{guolak}, where in our case the nonlinearity can be allowed to be discontinuous.
\end{remark}
\section*{Acknowledgements}
The authors would like to thank the anonymous Referee for the careful reading of the manuscript and
the constructive comments.
A. Calamai and G.~Infante are members of the Gruppo Nazionale per l'Analisi Matematica, la Probabilit\`a e le loro Applicazioni (GNAMPA) of the Istituto Nazionale di Alta Matematica (INdAM).
G.~Infante is a member of the UMI Group TAA  ``Approximation Theory and Applications''. 
J. Rodr\'iguez--L\'opez has been partially supported by the VIS Program of the University of Calabria, by Ministerio de Ciencia y Tecnología (Spain), AEI and Feder, grant PID2020-113275GB-I00, and by Xunta de Galicia, grant ED431C 2023/12.
This study was partly funded by: Research project of MIUR (Italian Ministry of Education, University and Research) Prin 2022 “Nonlinear differential problems with applications to real phenomena” (Grant Number: 2022ZXZTN2).

\textsc{Alessandro Calamai, 
Dipartimento di Ingegneria Civile, Edile e Architettura,
Universit\`{a} Politecnica delle Marche,
Via Brecce Bianche,
60131 Ancona, Italy}\\
\emph{Email address:} \texttt{a.calamai@univpm.it}

\textsc{Gennaro Infante, Dipartimento di Matematica e Informatica, Universit\`{a} della
Ca\-la\-bria, 87036 Arcavacata di Rende, Cosenza, Italy}\\
\emph{Email address:} \texttt{gennaro.infante@unical.it}

\textsc{Jorge Rodr\'iguez--L\'opez, CITMAga \& Departamento de Estat\'{\i}stica, An\'alise Matem\'atica e Optimizaci\'on, Universidade de Santiago de Compostela,  15782, Facultade de Matem\'aticas, Campus Vida, Santiago, Spain}\\
\emph{Email address:} \texttt{jorgerodriguez.lopez@usc.es}

\end{document}